\documentclass[12pt,reqno]{amsart}
\usepackage[utf8]{inputenc}
\usepackage[english]{babel}
\usepackage{comment}
\usepackage{amsmath,amssymb,amsfonts,amsbsy,amsthm,amscd,latexsym,graphicx}
\usepackage{empheq,textcomp}
\usepackage{xcolor}
\usepackage{hyperref}

\theoremstyle{definition}
\newtheorem{theorem}{Theorem} [section]
\newtheorem{corollary}[theorem]{Corollary}
\newtheorem{lemma}[theorem]{Lemma}

\newtheorem{remark}[theorem]{Remark}
\newtheorem{example}[theorem]{Example}

\numberwithin{equation}{section}

\AtEndDocument{
  \par
  \bigskip
  \begin{tabular}{@{}l@{}}%
    \text{Department of Mathematics, University of Oregon, Eugene, OR 97403-1222, USA}\\
    \textit{E-mail address}: \texttt{putingyu@uoregon.edu}
  \end{tabular}}

\newcommand{\C}{\mathbb{C}}
\newcommand{\Ac}{{\mathcal{A}}}

\newcommand{\N}{\mathbb{N}}
\newcommand{\R}{\mathbb{R}}

\newcommand{\Z}{\mathbb{Z}}

\newcommand{\Eq}{\, = \,}

\newcommand{\Le}{\, \le \,}

\newcommand{\Ge}{\, \geq \,}

\newcommand{\qeddef}{{\quad $\diamondsuit$}}

\newcommand{\bigabs}[1]{\bigl|\,#1\,\bigr|}
\newcommand{\Bigabs}[1]{\Bigl|\,#1\,\Bigr|}

\newcommand{\ip}[2]{\langle\,#1,#2\,\rangle}
\newcommand{\bigip}[2]{\bigl\langle \,#1, \, #2 \,\bigr\rangle}
\newcommand{\Bigip}[2]{\Bigl\langle \,#1, \, #2 \,\Bigr\rangle}

\newcommand{\norm}[1]{\|\,#1\,\|}

\newcommand{\Bignorm}[1]{\Bigl\|\,#1\,\Bigr\|}

\newcommand{\bigparen}[1]{\bigl(\,#1\,\bigr)}
\newcommand{\Bigparen}[1]{\Bigl(\,#1\,\Bigr)}

\newcommand{\set}[1]{\{#1\}}
\newcommand{\bigset}[1]{\bigl\{#1\bigr\}}
\newcommand{\Bigset}[1]{\Bigl\{#1\Bigr\}}

\newcommand{\clspan}{{\overline{\text{span}}}}

\newcommand{\inN}{_{n\in\N}}
\newcommand{\sumli}{\sum_{n=1}^\infty}

\newcommand{\Lm}{L^2(\mu)}
\newcommand{\itnwon}{\set{A^nx}_{n\geq0,\,x\in S}}
\newcommand{\itnwo}{\frac{A^nx}{\norm{A^nx}}}
\newcommand{\itn}{\bigset{\frac{A^nx}{\norm{A^nx}}}_{n\geq 0}}
\begin{document}
\title{Spectrum of normal operators that generate
certain scalable iterative systems}
\author{Pu-Ting Yu}
%\affil[1]{aaa}
\date{September 2025}
\begin{abstract}
Let $A\colon H\rightarrow H$ be a normal operator on an infinite-dimensional separable Hilbert space $H$ and let $S\subseteq H$ be a finite subset such that $\set{A^nx}_{n\geq 0,\,x\in S}$ can be rescaled to form a frame for $H$. 
That is, there exist some subsets $J_x\subseteq \N\cup\set{0}$ and  some set of nonzero scalars $(c_{n,x})_{n\in J_x,\,x\in S}$ such that $\set{c_{n,x}A^nx}_{n\in J_x,\,x\in S}$ forms a frame for $H.$ Assume that there exist some $\eta\in\N$ and $\delta>0$ such that for each infinite $J_x$ there is an increasing syndetic subsequence $(n^x_{k})_{k\in \N}\subseteq J_x$ satisfying $|c_{n^x_{k},x}|\norm{A^{i^x_{k}}x}\geq \delta$ for some non-negative integers $i^x_{k}$ with $|i^x_{k}- n^x_{k}|\leq \eta$ for all $k\in \N$. We prove that there exist finitely many numbers $(r_i)_{i=1}^N$ such that the continuous spectrum of $A$ is concentrated on arcs of a circle centered at origin with radius $r_i$. In particular, $A$ must be a diagonal operator if $S$ is a singleton.

As an application, we establish the conjecture proposed by Aldroubi et al.\ asserting that the iterative system $\set{\frac{A^nx}{\norm{A^nx}}}_{n\geq 0,\,x\in S}$ is never a frame for $H$, provided one of the following two conditions holds: (i) The continuous spectrum of $A$ contains more than $|S|-1$ points with distinct moduli; (ii) $S$ is a singleton and $A$ is not a diagonal operator. %Partial results regarding the case that $S$ is countably infinite are also presented.    
\end{abstract}
    
\pagestyle{plain}
\maketitle

\section{Introduction}
Let $H$ be an infinite-dimensional separable Hilbert space equipped with inner product $\ip{\cdot}{\cdot}$ and let $A\colon H\rightarrow H$ be a bounded linear operator. The \emph{iterative system} generated by $A$ and a countable subset $S\subseteq H$ is defined as  $\set{A^nx}_{n\geq0,\,x\in S}.$ The study of such systems was motivated by the problem of determining conditions under which time samples (obtained by iterating $A$) can be traded for spatial samples in an evolving system (see \cite{AG01}, \cite{BG13}, \cite{CRT06}, \cite{LM94} \cite{QS06}, and references therein). This problem was completely solved in the finite-dimensional setting by Aldroubi et al.\ in \cite{ACMT17}, whereas a solution in full generality remains unknown. Assuming $A$ is a self-adjoint diagonal operator, the question as to whether such an iterative system can exhibit certain approximation properties -- such as completeness, minimality or frame property -- was also investigated in \cite{ACMT17}. 
Here we say a sequence $\set{x_n}\inN\subseteq H$ is complete and minimal if its closed span is equal to $H$ and there exists a sequence $\set{y_n}\inN\subseteq H$ such that $\ip{x_m}{y_n}=\delta_{mn}.$ Finding conditions on $A$ and $S$ that ensures the associated iterative system forms a frame for $H$ is a more intractable problem and is still widely open. A sequence $\set{x_n}\inN\subseteq H$ is said to be a \emph{frame} for $H$ if there exist positive constants $B$ and $C$, called \emph{frame bounds}, such that 
\begin{equation}
\label{frame_ineq}
B\,\norm{x}^2\Le \sumli|\ip{x}{x_n}|^2\Le C\,\norm{x}^2\quad\text{for all }x\in H.
\end{equation}
We say that $\set{x_n}\inN$ is a \emph{Bessel sequence} if at least the upper inequality of Inequality (\ref{frame_ineq}) is satisfied. For other studies that investigate the approximation properties of iterative systems generated by generic bounded linear operators or some other specific linear operators, see \cite{ACMT17}, \cite{ACKM25}, \cite{VB24}, \cite{CHR18}, \cite{CH19} and \cite{CHPS24}. Building upon the work established in \cite{ACMT17}, Aldroubi et al.\ further studied this problem for general (bounded) normal operators in \cite{ACMCP17}. It was shown in this paper that every iterative system generated by a normal operator cannot be simultaneously complete and minimal in $H$, regardless of the choice of $S$. Assuming that $S$ is a singleton, a complete characterization of iterative system generated by normal operators that can be a frame for $H$ was proved in \cite{ACMCP17}. This result was later extended by Cabrelli et al.\ in \cite{CMPP20} to the case where $S$ is an arbitrary finite subset of $H$. However, the case that $S$ is an arbitrary countable subset of $H$ remains open. We refer to \cite{Ph17} and \cite{CMS21} for other research papers relevant to iterative systems generated by normal operators.

While the conditions under which $\itnwon$ forms a frame for $H$ when $S$ is finite and $A$ is normal are known, a characterization of 
such an iterative system that can be rescaled to form a frame for $H$ is still unknown. That is, under what conditions there exist some sequences of scalars $(c_{n,x})_{n\geq0\,,x\in S}$ such that $\set{c_{n,x}A^nx}_{n\geq0,\,x\in S}$ is frame for $H$. It was conjectured by Aldroubi in \cite{ACMCP17} that $\bigset{\frac{A^nx}{\norm{A^nx}}}_{n\geq 0,\,x\in S}$ can never be a frame for $H$ if $A$ is a normal operator. %Following the terminology in \cite{PY24}, this means that no iterative systems generated by normal operators, $\itnwon$, can be \emph{frame-normalizable}, regardless of the choice of $S$. 
%Here a sequence $\set{x_n}\inN\subseteq H$ is defined as frame-normalizable if the normalized sequence $\bigset{\frac{x_n}{\norm{x_n}}}\inN$ is a frame for $H.$ 
In the same paper where this conjecture was proposed,  the authors established it for self-adjoint operators. Assuming that $S$ is finite, the author proved in \cite{PY24} that this conjecture is true if $\itnwon$ is a frame for $H$. To the best of our knowledge, the literature currently contains only a limited number of results addressing this conjecture. It is noteworthy that there have been many  results concerning sequences in $H$ that can be rescaled to form a frame for $H$. For example, Kutyniok et al.\ studied sequences that can be rescaled to yied a \emph{Parseval frame} (i.e. a frame with bounds equal to $1$). Sequences that can normalized to form a Bessel sequence or a frame was investigated by the author in \cite{PY24}. Additionally, the question of determining under what conditions can one construct a frame by sampling a sequence that can be rescaled to form a frame were considered in \cite{FS19}, \cite{Bo24} and \cite{BY25}, where the results are crucial for the discretization problem for continuous frames.  

In this paper, we investigate spectrum of normal operators that generate iterative system that can be rescaled to yield a frame. Overall, our main results are summarized in the following two statements. Let $S$ be an arbitrary finite subset of $H$ and let $(c_{n,x})_{n\in J_x,\,x\in S}$ be a sequence of nonzero scalars such that $\set{c_{n,x}A^nx}_{n\in J_x,\,x\in S}$ is a frame for $H$.
\begin{enumerate}
    \item [\textup(a)]  Assume that there exist some $\eta\in\N$ and $\delta>0$ such that for each infinite $J_x$ there is an increasing syndetic subsequence $(n^x_{k})_{k\in \N}\subseteq J_x$ satisfying $|c_{n^x_{k},x}|\norm{A^{i^x_{k}}x}\geq \delta$ for some non-negative integers $i^x_{k}$ with $|i^x_{k}- n^x_{k}|\leq \eta$ for all $k\in \N$. Then there exist finitely many numbers $(r_i)_{i=1}^{N}$ with $N\leq |S|-1$ such that such that the continuous spectrum of $A$ is concentrated on arcs of a circle centered at origin with radius $r_i$.  Moreover, $A$ must be a diagonal operator if $S$ is a singleton.
\end{enumerate}
We then apply statement (a) to partially prove the conjecture proposed by Aldroubi et al.\ in \cite[Remark 6.3]{ACMCP17} as follows.
\begin{enumerate}
    \item [\textup(b)] 
    Assume that there exist more than $|S|-1$ complex numbers in the continuous spectrum of $A$ with distinct moduli. Then for any $\set{x_i}_{i=1}^N\subseteq H$ and any syndetic subsequences $\set{n_k^1}_{k\in \N},\dots,\set{n_k^N}_{k\in \N}$ of $\N\cup\set{0}$ , the iterative system $\Bigset{\frac{A^{n_k^i}x}{\norm{A^{n_k^i}x}}}_{k\in\N,\,1\leq i\leq N}$ is not a frame for $H$. In the case that $S$ is a singleton, $\Bigset{\frac{A^{n_k^1}x}{\norm{A^{n_k^1}x_1}}}_{k\in\N}$ is never a frame for $H$ if $A$ is not a diagonal operator.
\end{enumerate} 

The remainder of this paper is structured as follows. Section 2 presents necessary definitions and some known results that will be used throughout this paper. In Section 3, we prove our main theorems and conclude this paper by some illustrative examples.

\section{Preliminaries}
Throughout this paper, all normal operator is assumed to be bounded and we use $A$ to denote a normal operator from $H$ to $H$. Unless otherwise specified, $S$ means an arbitrary finite subset of $H.$ The notation $\sigma(A)$ stands for the \emph{spectrum of $A.$} A complex number $\lambda\in \sigma(A)$ is called a \emph{spectral value} of $A$. It is well known that a spectral value of $A$ is either an eigenvalue or belongs to the continuous spectrum of $A$. By \emph{continuous spectrum} we mean the set of all spectral values $\lambda$ for which $A-\lambda I$ is one-to-one and has dense range in $H$, but its range is not equal to $H.$
Since we are particularly interested in iterative systems that can rescaled to form a frame for $H$, we denote by $$\Ac\bigparen{A,S,J_x,(c_{n,x})}\coloneq\set{c_{n,x}A^nx}_{n\in J_x,\,x\in S}$$
, where $J_x\subseteq\N\cup\set{0}$ are increasing sequences and $(c_{n,x})_{n\in J_x}$ are sequences of scalars associated with $x\in S$.
By removing corresponding $n$ from $J_x$ if necessary, we assume that $c_{n,x}\neq 0$ for all $n\in J_x$ and all $x\in S.$ In the case where all $J_x$ coincide with some common set $J$, we will simply write $\Ac\bigparen{A,S,J_x,(c_{n,x})}.$ We say a sequence of integers $(c_n)\inN$ is \emph{syndetic} if $\sup_{n\in\N}|c_{n+1}-c_n|<\infty.$ 

Every normal operator admits a spectral decomposition. Following the notations in \cite{Co90}, we state the \emph{spectral theorem} (with multiplicity) for normal operators as follows. The proof of this theorem can be found in \cite[Section 10, Chapter IX]{Co90}. Given a non-negative regular Borel measure on $\C$ with compact support, we define the multiplication operator associated with $\mu$ by $$(N_\mu f)(z)=zf(z),\quad z\in \text{supp}(\mu).$$ 
\begin{theorem}[Spectral theorem]
\label{spectral_thm}Let $A\colon H\rightarrow H$ be a normal operator. Then there exist a sequence of mutually singular non-negative Borel measures $\set{\mu_j}_{j\in \N\cup \set{\infty}}$ that are supported in $\sigma(A)$ and an unitary operator $$U\colon H\rightarrow \bigparen{L^2(\mu_\infty)}^{\infty}\oplus L^2(\mu_1)\oplus \bigparen{L^2(\mu_2)}^{(2)}\oplus \cdots$$
 such that 
 $$UAU^{-1} \Eq N_{\mu_{\infty}}^{(\infty)}\oplus N_{\mu_{1}}\oplus N_{\mu_{2}}^{(2)}\oplus\cdots $$
 Moreover, if $M$ is another normal operator with corresponding Borel measure $\set{\nu_j}_{j\in \N\cup\set{\infty}}$, then $M$ is unitarily equivalent to $A$ if and only if $[\mu_j]=[\nu_j].$   
 \qeddef
\end{theorem}
For notational simplicity, for a given normal operator $A$, we denote by $V_A$ the associated Hilbert space in the spectral theorem, i.e., $$V_A=\bigparen{L^2(\mu_\infty)}^{\infty}\oplus L^2(\mu_1)\oplus \bigparen{L^2(\mu_2)}^{(2)}\oplus \cdots.$$ Moreover, since $\mu_j$ are mutually singular, we define the \emph{scalar spectral measure} $\mu_A$ associated with $A$ by 
$$\mu_A\coloneq \sum_{j=1}^\infty\mu_j.$$
To avoid ambiguity when switching between the norms associated with $H$ and $V_A$, we use the notations $\norm{\cdot}$ and $\norm{\cdot}_{L^2(\mu)}$ to denote the norms associated with $H$ and $V_A$, respectively.
%\textcolor{red}{If $\Bigset{\frac{A^{n_k}x}{\norm{A^{n_k}x}}}_{k\in\N}$ is a frame for $H$ for some syndetic sequence $\set{n_k}_{k\in \N}$, then $\itn$ is a frame for $H$.}

Let $\mathbb{D}=\set{z\,|\,|z|<1}$ be the open unit disk in $\C$ and let $\text{Hol}(\mathbb{D})$ be the set of all holomorphic functions defined on $\mathbb{D}$. The Hardy space $H^2(\mathbb{D})$ is the Hilbert space $$H^2(\mathbb{D})=\set{\phi\in \text{Hol}(\mathbb{D})\,|\,\phi=\sum_{n=0}^\infty c_nz^n\text{ for some }(c_n)_{n=0}^\infty\in \ell^2(\N\cup\set{0})},$$
equipped with the inner product $\ip{\phi}{\psi}=\sum_{n=0}^\infty c_n\overline{d_n}$, provided that $\phi=\sum_{n=0}^\infty c_nz^n$ and $\psi =\sum_{n=0}^\infty d_nz^n$. Given a sequence of $\Lambda=(\lambda_n)_{n\in\N}\subseteq \mathbb{D}$, the associated \emph{evaluation operator} with $(\lambda_n)_{n\in\N}$, 
$T_{\Lambda}\colon D(T_\Lambda)\subseteq H^2(\mathbb{D})\rightarrow \ell^2(\N),$ is defined by $$T_{\Lambda}\phi\coloneq \Bigparen{\sqrt{1-|\lambda_n|^2}\phi(\lambda_n)}\inN$$, where $D(T_\Lambda)$ denotes the domain of $T_{\Lambda}.$ If $(x_n)\inN\subseteq H$ is a Bessel sequence in $H$, then the corresponding \emph{analysis operator} is the bounded linear operator $$C_{(x_n)}\colon H\rightarrow \ell^2(\N), \quad C_{(x_n)}(x)=\bigparen{\ip{x}{x_n}}\inN.$$
 A sequence of scalars $(\lambda_n)\inN\subseteq \mathbb{D}$ is said to be \emph{uniformly separated (or satisfies the \emph{Carleson condition})} if $$\inf\limits_{n} \prod_{k\neq n} \dfrac{|\lambda_n-\lambda_k|}{\bigabs{1-\overline{\lambda_n}\lambda_k}}\geq \delta>0.$$ 
We now summarize some known results below regarding the iterative system $\bigset{\itnwo}_{n\in\N,\,x\in S}$ that will be used throughout this paper.

\begin{theorem}
\label{required_results}
Let $A\colon H\rightarrow H$ be a normal operator.
     \begin{enumerate}
    \setlength\itemsep{0.5em}
    \item[\textup{(a)}] (\cite[Proposition 4.1]{ACMCP17}) The iterative system $\itnwon$ is never complete and minimal in $H$ for any countable subset $S$ of $H$.
     \item[\textup{(b)}] (\cite[Theorem 2.3]{ACMT17}) Let $S=\set{x_i}_{i=1}^N$ be a finite subset of $H$. The iterative system $\set{A^nx_i}_{n\geq 0,\,1\leq i\leq N}$ is a frame for $H$ if and only if the following conditions hold:
     \begin{enumerate}
        \setlength\itemsep{0.5em}
         \item[\textup{(i)}] $A=\sum_{n\in\N}\lambda_nP_n$, where $P_n$ are orthogonal projections onto the eigenspaces associated with $\lambda_n.$ Moreover, $\text{rank}(P_n)\leq |S|$ for all $n\in \N,$
         \item[\textup{(ii)}] $(\lambda_n)\inN$ is a union of $|S|$ uniformly separated sequences in $\mathbb{D},$
         \item[\textup{(iii)}]There exist some positive constants $\alpha$ and $\beta$ such that for all $n\in\N$ and all $y\in P_nH$ we have 
         $$\alpha (1-|\lambda_n|^2)\norm{y}^2 \Le \sum _{i=1}^N\,\bigabs{\ip{y}{P_nx_i}}^2\Le \beta (1-|\lambda_n|^2)\norm{y}^2,$$
         \item[\textup{(iv)}] $\bigparen{\text{Range}(T_{\Lambda})}^{N}+\text{Ker}(C^*_{E})\Eq \ell^2(N\times \N)$, where for $x\in H$, $$C_E(x)\Eq \Bigparen{\bigparen{\ip{x}{(1-|\lambda_n|^2)^{-\frac{1}{2}}P_nx_1}}\inN,\dots,\bigparen{\ip{x}{(1-|\lambda_n|^2)^{-\frac{1}{2}}P_nx_N}}\inN}.$$ 
     \end{enumerate}
     \item[\textup{(c)}](\cite[Theorem 5.7]{PY24}) If $\itnwon$ is a frame for $H$, then for any finite subset $S\subseteq H$ the normalized iterative system $\bigset{\itnwo}_{n\geq0,\,x\in S}$  is not a frame for $H.$ 
    \end{enumerate}
    We remark that in the case that $S=\set{x_1}$ condition (ii) in Theorem \ref{required_results} is equivalent to $$\alpha\leq\dfrac{\norm{P_nx_1}^2}{(1-|\lambda_n|^2)}\leq \beta, \quad \text{for all }n\in\N.$$
    Moreover, statement (iv) becomes redundant in the case where $S$ is a singleton, since condition (ii) then implies condition (iv) (see \cite[Remark 2.4]{ACMT17} for further comparisons between the cases where $S$ is a singleton and where $S$ is an arbitrary finite subset of $H$). 
\section{Main Results}
In summary, our goal is to characterize the spectrum of normal operators that generate certain iterative systems which, after appropriate rescaling, form frames for $H.$ Specifically, iterative systems that satisfy the following conditions.
\begin{enumerate}
\setlength\itemsep{0.5em}
    \item [\textup{(a)}] There exists some subsets $J_x\subseteq \N\cup\set{0}$ and some sequences of scalars $(c_{n,x})_{n\in J_x,\,x\in S}$ such that $\Ac\bigparen{A,S,J_x,(c_{n,x})}$ is a frame for $H.$
     \item [\textup{(b)}] There exist $\eta\in\N$ and $\delta>0$ such that for each infinite $J_x$ there is an increasing syndetic subsequence $(n^x_{k})_{k\in \N}\subseteq J_x$ satisfying $|c_{n^x_{k},x}|\norm{A^{i^x_{k}}x}\geq \delta$ for some non-negative integers $i^x_{k}$ with $|i^x_{k}- n^x_{k}|\leq \eta$ for all $k\in \N$.
\end{enumerate}
We would like point out that sequences of scalars that satisfy condition (b) are quite common. For instance, when $J=\N\cup\set{0}$ and $c_{n}=\norm{A^nx}^{-1}$, we obtain the normalized iterative system $\Ac\bigparen{A,S,J_x,(c_{n,x})}=\bigset{\frac{A^nx}{\norm{A^nx}}}_{n\geq0}.$
Let $(i_n)_{n\geq0}\subseteq\N\cup\set{0}$ be any sequence such that $\sup_{n\in\N}|i_n-n|<\infty.$ As we will see in Lemma \ref{syndetic_frame_property}, the sequence $\bigparen{\norm{A^{i_n}x}^{-1}}_{n\in\N}$ also satisfies condition $(b)$, which gives rise to the system $\bigset{\frac{A^nx}{\norm{A^{i_n}x}}}_{n\geq0}$. 

%beginning with a series of lemmas. 

We begin with a series of lemmas. By Theorem \ref{spectral_thm}, normal operators behave like multiplication operators on Hilbert spaces. Consequently, we obtain the first lemma as follows.
Note that for any normal operator $A$ we have $\text{ker}(A)=\text{ker}(A^n)$ for any $n\geq 2.$ Consequnetly, for any $y\in H$ we have that $A^ny=0$ for some $n\geq 2$ if and only if $Ay=0$. 
 \begin{lemma}
\label{fraction_conv_operator_norm}
     Let $A\colon H\rightarrow H$ be a normal operator. Then for any nonzero $y\in H$ such that $Ay\neq 0$ 
     we have $$\frac{\norm{A^{n+1}y}}{\norm{A^ny}}\nearrow\sup_{z\in \text{supp}(Uy)}|z|\quad \text{as }n\rightarrow\infty.  $$
\end{lemma}
\begin{proof}

Let $y\in H\setminus\text{ker}(A)$.
   %Since $\set{A^nx}_{n\geq 0}$ is complete in $H$, we see that $A^n x\neq 0$ for any $n\in \N.$ 
   By the Cauchy-Bunyakovski-Schwarz inequality, we have 
   \begin{equation}
   \label{limit_operator_norm}
   \norm{Ay}^2\Eq \ip{y}{A^*Ay}\Le \norm{A^*Ay}\norm{y}\Eq \norm{A^2{y}}\norm{y},
    \end{equation}
    %Then for any $y\in H$, we obtain $$\frac{\norm{Ay}}{\norm{y}}\Le \frac{\norm{A^2y}}{\norm{Ay}}$$
   which show that $\frac{\norm{Ay}}{\norm{y}}\leq \frac{\norm{A^2y}}{\norm{Ay}}$. By considering $y'=Ay$ in Equation (\ref{limit_operator_norm}), we obtain that $\frac{\norm{A^2y}}{\norm{Ay}}\leq \frac{\norm{A^3y}}{\norm{A^2y}}$. Inductively, we see that $\bigset{\frac{\norm{A^{n+1}y}}{\norm{A^ny}}}_{n\geq0}$ is an increasing sequence. Therefore, $\lim\limits_{n\rightarrow \infty }\frac{\norm{A^{n+1}y}}{\norm{A^ny}}$ exists. %Since $\set{A^nx}_{n\geq 0}$ is complete in $H$, it follows that $\bigset{\frac{\norm{A^{n+1}x}}{\norm{A^nx}}}_{n\geq0}$ is an increasing sequence.
%such that for any $n\in \N$ there exists $n_\delta$ such that $\frac{\norm{A^{n_\delta+1}y}}{\norm{A^{n_\delta}y}}\leq C-2\delta.$

   Next, let $C=\sup_{z\in \text{supp}(Uy)}|z|$ and let $0<\delta<C$. By assumption, the set $$E=\bigset{z\in \text{supp}(Uy)\,\big|\,|z|\geq C-\delta}$$ has a positive $\mu$-measure. Therefore, for each $n\geq 0$ we have 
   $$\frac{\norm{A^{n+1}y}^2}{\norm{A^ny}^2}=\frac{\norm{N_\mu^{n+1} \bigparen{(Uy)\cdot\chi_E}}^2_{\Lm}+\norm{N_\mu^{n+1} \bigparen{(Uy)\cdot\chi_{E^c}}}^2_{\Lm}}{\norm{N_\mu^{n} \bigparen{(Uy)\cdot\chi_E}}^2_{\Lm}+\norm{N_\mu^{n} \bigparen{(Uy)\cdot\chi_{E^c}}}^2_{\Lm}}.$$
   By dominate convergence theorem, we see that $$\Bignorm{\frac{N_\mu^{n+1} \bigparen{(Uy)\cdot\chi_{E^c}}}{(\norm{A}-\delta)^n}}^2\rightarrow 0 \quad \text{and} \quad \Bignorm{\frac{N_\mu^{n} \bigparen{(Uy)\cdot\chi_{E^c}}}{(\norm{A}-\delta)^n})}^2\rightarrow 0\quad \text{as } n\rightarrow \infty.$$
   It then follows that 

   $$C^2\Ge\lim\limits_{n\rightarrow \infty }\frac{\norm{A^{n+1}y}^2}{\norm{A^ny}^2} \Ge \liminf\limits_{n\rightarrow \infty } \frac{\norm{N_\mu^{n+1} \bigparen{(Uy)\cdot\chi_E}}^2}{\norm{N_\mu^{n} \bigparen{(Uy)\cdot\chi_E}}^2}\Ge (C-\delta)^2.$$
 Since $\delta$ is arbitrary, we see that  $\lim\limits_{n\rightarrow \infty }\frac{\norm{A^{n+1}y}}{\norm{A^ny}}\Eq C.$
\end{proof}

\begin{lemma}
\label{complete_full_measure}
    Let $A\colon H\rightarrow H$ be a normal operator and let $S$ be a countable subset of $H$. Assume that $\set{A^{n_k}x}_{x\in S,\,k\in \N}$ is complete in $H$ for some sequence $\set{n_k}_{k\in \N}\subseteq \N\cup\set{0}$. Then $\cup_{x\in S}\, \text{supp}(Ux)=\text{supp}(\mu)$.
\end{lemma}
\begin{proof}
    Suppose to the contrary that $\cup_{x\in S}\, \text{supp}(Ux)\neq \text{supp}(\mu)$. Let $E\subseteq \text{supp}(\mu)$ be a set of positive $\mu$-measure such that $(Ux)|_{E}=0$ $\mu$-almost everywhere for all $x\in S.$ Then it follows that for each $k\in \N$ and each $x\in S$ 
    $$\ip{U^{-1}(\chi_{E})}{A^{n_k}x}\Eq \ip{\chi_E}{UA^{n_k}x}=\ip{\chi_E}{N_{\mu}^{n_k}(Ux)}=0,$$
    which is a contradiction.
\end{proof}

The next lemma proves a sufficient condition for an iterative system to be a frame for $H$.  
\begin{lemma}
\label{syndetic_frame_property}
 Let $A\colon H\rightarrow H$ be a normal operator such that $\mathcal{A}\bigparen{A,S,J_x,(c_{n,x})}$ is a frame for $H$ for some finite subset $S$. 
Assume that there exist some $\eta\in\N$ and $\delta>0$ such that for each infinite $J_x$ there is an increasing syndetic subsequence $(n^x_{k})_{k\in \N}\subseteq J_x$ satisfying $|c_{n^x_{k},x}|\norm{A^{i^x_{k}}x}\geq \delta$ for some non-negative integers $i^x_{k}$ with $|i^x_{k}- n^x_{k}|\leq \eta$ for all $k\in \N$.
 Then $\bigset{\frac{A^nx}{\norm{A^nx}}}_{n\geq 1,\,x\in S\setminus \text{ker}(A)}$ is a frame for $H$.
\end{lemma}
\begin{proof}
For notational convenience, we denote $S\setminus \text{ker}(A)$ by $S_1$.
Let $B$ be an upper frame bound of $\mathcal{A}(A,S,J_x,(c_{n,x}))$ and let $$L=\sup\limits_{\substack{k\in \N,\,x\in S_1}}|n^x_{k+1}-n^x_{k}|.$$ 
%We first consider the case that $S_1=S.$ 
Without loss of generality, we may assume that $n^x_{1}\geq \eta$ for all $x\in S_1$. We claim that there exists some constant $m'>0$ such that $ \inf\limits_{\substack{k\in\N,\,x\in S_1}}\frac{|c_{n^x_k,x}|}{|c_{n^x_{k+j},x}|}\geq m'$ for all $1\leq j\leq L$.
 Since $|c_{n^x_k,x}|\norm{A^{n^x_k}x}\leq B$ for all $k\in \N$ and all $x\in S$, for each $x\in S_1$ we have  
 %By Lemma \ref{fraction_conv_operator_norm}, $\bigparen{\frac{\norm{A^{n+1}y}}{\norm{A^{n}y}}}_{n\geq 0}$ is an increasing sequence for any $y\in H$. It then follows that for any $1\leq j\leq L$ 
\begin{equation}
    \frac{|c_{n^x_{k},x}|}{|c_{n^x_{k+j},x}|}\Ge \frac{\delta}{B} \frac{\norm{A^{n^x_{k+j}}x}}{\norm{A^{i^x_{k}}x}}\Eq \frac{\delta}{B} \frac{\norm{A^{n^x_{k+j}}x}}{\norm{A^{n^x_{k}}x}} \frac{\norm{A^{n^x_{k}}x}}{\norm{A^{i^x_{k}}x}}.
\end{equation}
We first evaluate  $\frac{\norm{A^{n^x_{k+j}}x}}{\norm{A^{n^x_{k}}x}}$. Using telescoping and Lemma \ref{fraction_conv_operator_norm}, we obtain
$$\frac{\norm{A^{n^x_{k+j}}x}}{\norm{A^{n^x_{k}}x}}\Eq \frac{\delta}{B}\frac{\norm{A^{n^x_{k+j}}x}}{\norm{A^{n^x_{k+j}-1}x}}\cdots
    \frac{\norm{A^{n^x_{k}+1}x}}{\norm{A^{n^x_{k}}x}}  \Ge \frac{\delta}{B}\Bigparen{\frac{\norm{Ax}}{\norm{x}}}^{n^x_{k+j}-n^x_k}.$$
Note that $n^x_{k+j}-n^x_k\leq L^2$ for all $x\in S$, all $k\in \N$ and all $1\leq j\leq L$. On the other hand, using the same telescoping argument, we see that if $i_k^x< n_k^x$, then $\frac{\norm{A^{n^x_{k}}x}}{\norm{A^{i^x_{k}}x}}\Ge \bigparen{\frac{\norm{Ax}}{\norm{x}}}^{n^x_{k}-i^x_k}.$
For the case that $i_k^x>n_k^x$, we have 
$$\frac{\norm{A^{n^x_{k}}x}}{\norm{A^{i^x_{k}}x}}=\frac{\norm{A^{n^x_{k}}x}}{\norm{A^{n^x_{k}+1}x}}\cdots \frac{\norm{A^{i^x_{k}-1}x}}{\norm{A^{i^x_{k}}x}}\Ge \norm{A}^{-(i_k^x-n_k^x)}.$$
We then set $m'=\frac{\delta}{B}\min\limits_{k\in\N,\,x\in S_1}\Bigparen{ \bigparen{\frac{\norm{Ax}}{\norm{x}}}^{n^x_{k+j}-n^x_k}\bigparen{\frac{\norm{Ax}}{\norm{x}}}^{n^x_{k}-i^x_k},\,\bigparen{\frac{\norm{Ax}}{\norm{x}}}^{n^x_{k+j}-n^x_k}\norm{A}^{n^x_{k}-i^x_k}}$.
%\Ge \frac{\delta}{B}\frac{\norm{A^{n_{k+j}}x}}{\norm{A^{n_{k+j}-1}x}}\cdots
%    \frac{\norm{A^{n_{k}+1}x}}{\norm{A^{n_{k}}x}}\Ge \frac{\delta}{B} \Bigparen{\frac{\norm{Ax}}{\norm{x}}}^{n_{k+j}-n_k}.
%On the other hand, we have  $$\frac{|c_{n_{k+1}}|}{|c_{n_{k}}|}\Ge \frac{\delta}{B}\frac{\norm{A^{n_{k+1}}x}}{\norm{A^{n_{k}}x}}\Ge \frac{\delta}{B\norm{A}^L}.$$

%Also, we may assume $A^nx\neq 0$ for all $n\in\N$ and all $x\in S$.
Next, without loss of generality, for each $x\in S_1$ we assume that $I_x=\set{n^x_k\,|\,k\in \N}$ is the maximal subset of $J_x$ such that $|c_{n^x_k,x}|\norm{A^{n^x_k}x}\geq \delta$ for all $k\in \N.$ We show that there exists some positive constant $m''$ such that 
\begin{align}
    \begin{split}
   \sum_{j=1}^L \sum_{x\in S_1}\sum_{k=1}^\infty \,\bigabs{\ip{y}{c_{n^x_k,x}A^{n^x_k+j}x}}^2 \Ge m'' \sum_{x\in S_1}\sum_{n\in J_x} \,\bigabs{\ip{y}{c_{n,x}A^{n}x}}^2
     \end{split}
    \end{align}
Fix $x\in S_1$ and fix $1\leq j\leq L$. For those $k$ such that $n_k^x+j\in I_x$ , we have that for any $y\in H$  
\begin{align*}
%\label{lower_frame_eq1}
\begin{split}
\bigabs{\ip{y}{c_{n^x_k,x}A^{n^x_{k}+j}x}}^2&\Eq \bigabs{\ip{y}{c_{n^x_{k}+j,x}A^{n^x_{k}+j}x}}^2\dfrac{|c_{n^x_{k},x}|^2}{|c_{n^x_{k}+j,x}|^2}
\Ge (m')^2\bigabs{\ip{y}{c_{n^x_{k,x}+j}A^{n^x_{k}+j}x}}^2.
\end{split}
\end{align*}
For those $k$ with $n^x_k+j\in J_x\setminus I_x$, likewise, we have

\begin{align}
\label{lower_frame_eq2}
\begin{split}
\bigabs{\ip{y}{c_{n^x_k,x}A^{n^x_{k}+j}x}}^2 &\Eq \bigabs{\ip{y}{c_{n^x_{k}+j,x}A^{n^x_{k+j}}x}}^2\dfrac{|c_{n^x_{k},x}|^2}{|c_{n^x_{k}+j,x}|^2}\\&\Eq \dfrac{|c_{n^x_{k},x}|^2\norm{A^{i^x_k}x}^2}{|c_{n^x_{k}+j,x}|^2\norm{A^{n^x_k+j}x}^2} \Bigparen{\frac{\norm{A^{n^x_k+j}x}}{\norm{A^{i^x_k}x}}}^{2}\bigabs{\ip{y}{c_{n^x_{k}+j,x}A^{n^x_{k}+j}x}}^2\\
&\Ge \Bigparen{\frac{\norm{A^{n^x_k+j}x}}{\norm{A^{i^x_k}x}}}^{2}\bigabs{\ip{y}{c_{n^x_{k}+j,x}A^{n^x_{k}+j}x}}^2.
\end{split}
\end{align}

Here we used the fact that $|c_{n^x_{k}+j,x}|\norm{A^{n^x_k+j}x}<\delta$ to obtain the third inequality from the second inequality. Arguing in the same way as when deriving $m'$, we see that $\Bigparen{\frac{\norm{A^{n^x_k+j}x}}{\norm{A^{i^x_k}x}}}^{2}$ is bounded below by some positive number. Since for each $x\in S_1$ the sequence $(n_k^x)_{k\in\N}$ is syndetic, it follows that 
$\cup_{x\in S_1}J_x\subseteq \cup_{x\in S_1}\cup_{j=0}^L(I_x+j)$ for all $x\in S_1.$ 
Consequently, there exists some positive constant $m''$ such that 
\begin{align}
\label{lower_frame_estimate}
    \begin{split}
    \sum_{j=1}^L\sum_{x\in S} \sum_{k=1}^\infty\,\bigabs{\ip{y}{c_{n^x_k,x}A^{n_k^x+j}x}}^2 \Ge m'' \sum_{x\in S_1}\sum_{n\in J_x} \,\bigabs{\ip{y}{c_{n,x}A^{n}x}}^2
     \end{split}
    \end{align}
Let $n_0=\inf\limits_{k\in \N,\,x\in S_1} n_{k}^x$. Note that for each $x\in S_1$ and each $n\geq n_0$ there exist at most $L$ $n^x_k$ such that $n^x_k+j=n$ for some $1\leq j\leq L$. Then for each $n\geq n_0$ and $x\in S_1$ we define $$d_{n,x}=|c_{n^x_{k_1},x}|+|c_{n^x_{k_2},x}|+\cdots,$$
where $n^x_{k_i}+j_i=n$ for some $1\leq j_i\leq L.$  
By Equation (\ref{lower_frame_estimate}), we see that 
$$\sum_{x\in S_1}\sum_{n=n_0}^\infty \,\bigabs{\ip{y}{d_{n,x}A^nx}}^2 \Ge m'' \sum_{x\in S_1}\sum_{n\in J_x} \,\bigabs{\ip{y}{c_{n,x}A^{n}x}}^2.$$
Thus, $\mathcal{A}(A,S\setminus S_1,J_x,(c_{n,x}))\cup \bigset{d_{n,x}A^nx}_{n\geq n_0,\,x\in S}$ satisfies a lower frame condition in $H$. Moreover, by using telescoping if necessary, we see that there exists some $m>0$ such that $|d_{n,x}|\norm{A^nx}\geq m$ for all $n\geq n_0$ and all $x\in S_1.$ For notational simplicity, we denote $\mathcal{A}(A,S\setminus S_1,J_x,(c_{n,x}))$ by $\Ac_{S\setminus S_1}.$
By the boundedness of $A$, it follows that  $\Ac_{S\setminus S_1}\cup \bigset{d_{n,x}A^nx}_{n\geq n_0,\,x\in S}$ is a Bessel sequence in $H$, and hence a frame for $H$. 

It remains to show the removal of $\Ac_{S\setminus S_1}$ from $\Ac_{S\setminus S_1} \cup\bigset{d_{n,x}A^nx}_{n\geq n_0,\,x\in S}$ does not break the frame property. Note that $\Ac_{S\setminus S_1}=\set{c_{0,x}x}_{x\in S\setminus S_1}.$ By adding finitely many elements if necessary, we may assume that $n_0=1.$ It is known that the removal of an element from a frame either leaves a frame for $H$ or an incomplete set in $H$ (for example, see \cite[Theorem 8.22]{Hei11}). Suppose that the removal of some $x\in S\setminus S_1$ leaves an incomplete set. This assumption implies that $\Ac_{S\setminus S_1}\cup \bigset{d_{n,x}A^nx}_{n\geq n_0,\,x\in S}$ is complete and minimal in $H$, and hence $\set{x}_{x\in S\setminus S_1}\cup \set{A^nx}_{n\geq 1,\,x\in S_1}$ is complete and minimal in $H$. However, this would contradict Theorem \ref{required_results} (a). Consequently, the removal of any $x\in S\setminus S_1$ leaves a frame for $H$. Repeating this process at most finitely many times, we see that $\bigset{d_{n,x}A^nx}_{n\geq 1,\,x\in S_1}$ is a frame for $H$, and hence $(\norm{d_{n,x}A^nx})_{n\geq 1,\,x\in S_1}$ is bounded above by some constant. Since $\norm{d_{n,x}A^nx}=\frac{d_{n,x}}{\norm{A^nx}^{-1}}$ is bounded above and below by some positive constants, it follows that  $\bigset{\frac{A^nx}{\norm{A^nx}}}_{n\geq 1,\,x\in S_1}$ is a frame for $H$. \qedhere

%by adding finitely many elements into $\Ac_{S\setminus S_1}\cup \bigset{d_{n,x}A^nx}_{n\geq n_0,\,x\in S}$ if necessary, we may assume that 
%$J_x=\set{0,1,2,...,N}$ for some $N\in \N$ for all $x\in S\setminus S_1$ and assume that $n_0=0.$ 

\end{proof}
   
\end{theorem}
\begin{remark} 
 A special case of Lemma \ref{syndetic_frame_property} implies that $\mathcal{A}\bigparen{A,S,\N,(\norm{A^nx}^{-1})}$ is a frame for $H$ if there exist syndetic sequences $J_x$ such that $\mathcal{A}\bigparen{A,S,J_x,(\norm{A^nx}^{-1})}$ is a frame for $H.$ However, it is still unclear whether the frame property of $\mathcal{A}\bigparen{A,S,\N,(\norm{A^nx}^{-1})}$ would necessarily imply the frame property of $\mathcal{A}\bigparen{A,S,J_x,(\norm{A^nx}^{-1})}$ for any syndetic sequences $J_x.$ Nevertheless, using the fact that $\frac{A^nx}{\norm{A^nx}}$ all have the same norm together with Bownik's selector result (\cite[Theorem 5.3]{Bo24}), one can show that for any $\eta\in \N$ there exist syndetic subsequences $J_x$ satisfying $\sup\limits_{n_k\in J_x,x\in S}|n_{k+1}-n_k|\geq \eta$ such that $\mathcal{A}\bigparen{A,S,J_x,(\norm{A^nx}^{-1})}$ is a frame for $H.$ Since this is not the main interest of the present paper, we only provide a brief outline the construction as follows. Let $\alpha\in \N$ be chosen so that $2^\alpha$ is so large and let $L$ be such that $\frac{A^nx}{\norm{A^nx}L}$ is small enough. We then partition $\N$ into disjoint blocks $$\N=\cup_{k=0}^\infty J_k, \quad J_k=\set{n}_{n=k2^\alpha+1}^{(k+1)2^\alpha}.$$ For each $k\geq0$, we pair any two numbers in $J_k$ and applying \cite[Theorem 5.3]{Bo24}. Repeating this procedure $\alpha$ times yields the desired frame for $H$.   \qeddef
\end{remark}

The final lemma needed establishes the invertibility of $A$ under our primary assumptions.
\begin{lemma}
\label{invertible_operatorA}
 Let $A\colon H\rightarrow H$ be a normal operator such that $\mathcal{A}\bigparen{A,S,J_x,(c_{n,x})}$ is a frame for $H$ for some finite subset $S$. 
Assume that there exist some $\eta\in\N$ and $\delta>0$ such that for each infinite $J_x$ there is an increasing syndetic subsequence $(n^x_{k})_{k\in \N}\subseteq J_x$ satisfying $|c_{n^x_{k},x}|\norm{A^{i^x_{k}}x}\geq \delta$ for some non-negative integers $i^x_{k}$ with $|i^x_{k}- n^x_{k}|\leq \eta$ for all $k\in \N$. Then $A$ is invertible.
\end{lemma}
\begin{proof}
By Lemma \ref{syndetic_frame_property}, we may assume that $S\cap\text{ker}(A)=\emptyset.$
 Also, by Lemma \ref{syndetic_frame_property}, we have that $\mathcal{A}\bigparen{A,S,\N,(\norm{A^nx}^{-1})}$ is a frame for $H$. Since $A$ is normal, it is equivalent to show that $\text{Range}(A)=H.$ Let $\set{a_{n,x}^*}_{n\geq 1,\,x\in S}$ be the canonical dual frame associated with $\itn$. Then for any $y\in H$ we have 
\begin{align*}
\begin{split}   
y\Eq\sum_{x\in S} \sumli \ip{y}{a_{n,x}^*}\frac{A^nx}{\norm{A^nx}}\Eq A\Bigparen{ \sum_{x\in S} \sumli\ip{y}{a_{n,x}^*}\frac{A^{n-1}x}{\norm{A^nx}}}.
\end{split} 
\end{align*}
It remains to show that $ \sum_{x\in S}\sumli\ip{y}{a_{n,x}^*}\frac{A^{n-1}x}{\norm{A^nx}}$ converges in the norm of $H.$ Note by Lemma \ref{fraction_conv_operator_norm}, we have $\frac{\norm{A^{n}x}}{\norm{A^{n+1}x}}\leq \frac{\norm{x}}{\norm{Ax}}$ for all $n\geq 0$. It then follows that for any $y\in H$
\begin{align}
    \begin{split}
\sum_{x\in S}\sum_{n=1}^\infty \,\Bigabs{\Bigip{y}{\frac{A^{n-1}x}{\norm{A^nx}}}}^2&\Eq \sum_{x\in S}\sum_{n=1}^\infty \,\Bigabs{\Bigip{y}{\frac{A^{n-1}x}{\norm{A^{n-1}x}}}\frac{\norm{A^{n-1}x}}{\norm{A^{n}x}}}^2\\
&\Le \frac{\norm{x}}{\norm{Ax}}\sum_{x\in S}\sum_{n=1}^\infty \,\Bigabs{\Bigip{y}{\frac{A^{n-1}x}{\norm{A^{n-1}x}}}}^2.
\end{split}
\end{align}
Thus,
$\bigset{\frac{A^{n-1}x}{\norm{A^nx}}}_{n\geq 1}$ is a Bessel sequence.    
%Next, let $C>0$ be such that $\norm{CA}<\frac{1}{2}$ and let $A'=CA.$ Note that for each $n\geq 1$ we have $\frac{(A')^nx}{\norm{(A')^nx}}=\frac{A^nx}{\norm{A^nx}}$. Consequently, $\Bigset{\frac{(A')^nx}{\norm{(A')^nx}}}_{n\geq 1}$ is a still a frame for $H.$ Let $\set{a_n^*}_{n\in \N}$ b
\end{proof}
Given $r>0$, by $B_r(0)$ we mean the open ball in $\C$ centered at the origin with radisu $r.$
\begin{theorem}
\label{countable_spec_inside}
  Let $A\colon H\rightarrow H$ be a normal operator such that $\mathcal{A}\bigparen{A,S,J_x,(c_{n,x})}$ is a frame for $H$ for some finite subset $S$. 
Assume that there exist some $\eta\in\N$ and $\delta>0$ such that for each infinite $J_x$ there is an increasing syndetic subsequence $(n^x_{k})_{k\in \N}\subseteq J_x$ satisfying $|c_{n^x_{k},x}|\norm{A^{i^x_{k}}x}\geq \delta$ for some non-negative integers $i^x_{k}$ with $|i^x_{k}- n^x_{k}|\leq \eta$ for all $k\in \N$.  Then there exist some positive integer $N\leq |S\setminus\text{ker}(A)|$ and finitely radii $(r_i)_{i=1}^N\subseteq (0,\norm{A}]$ such that the following statements hold:
  \begin{enumerate}
  \setlength\itemsep{0.5em}
      \item [\textup{(a)}] The continuous spectrum of $A$ lies entirely on arcs of circles centered at origin with radii $r_i$, where $r_i$ are positive numbers for which $$S_i=\bigset{x\in S\,|\sup_{z\in \text{supp}(Ux)}|z|\Eq r_i\,}\neq\emptyset$$
       \item [\textup{(b)}] For any sufficiently small $\delta>0$, there exist at most finitely many spectral values of $A$ contained in $\overline{B_{r_1-\delta}}(0)\cup \Bigparen{\cup_{i=2}^N\bigparen{\overline{B_{r_i-\delta}}(0)\setminus B_{r_{i-1}+\delta}(0)}}$
  \end{enumerate}
  
\end{theorem}
\begin{proof}
By Lemma \ref{syndetic_frame_property}, we may assume that $S\cap\text{ker}(A)=\emptyset$, $J_x=\N$ and $c_{n,x}=\norm{A^nx}^{-1}$ for all $n\in \N$ and $x\in S.$
Let $r_0=0$ and let $r_1<r_2<\cdots <r_N$ be finitely many numbers for which $S_i=\bigset{x\in S\,|\sup_{z\in \text{supp}(Ux)}|z|\Eq r_i\,}\neq\emptyset\text{ for }1\leq i\leq N.$ 
Note that by Lemma \ref{complete_full_measure} and the fact that $A$ is a normal operator, we have $r_N=\norm{A}$.
Let $R_0=\emptyset$ and for $1\leq i\leq N$ we define $R_i=\bigset{x\in S\,|\sup_{z\in \text{supp}(Ux)}|z|\geq r_i\,}.$ %Note that $R_1\supseteq R_2\supseteq\cdots\supseteq R_N.$ 

Next, we show that for each $1\leq i\leq N$, the set $\sigma(A)\cap \bigparen{B_{r_{i}}(0)\setminus \overline{B_{r_{i-1}}}(0)}$ contains at most countably many discrete points. We first claim that for any $K\in \N$ (in fact, any $K\geq0$) the normalized iterative system $\bigset{\itnwo}_{n\geq K,\,x\in S}$ is also a frame for $H$. It is known that bounded bijective operators preserves the completeness of a complete set. Consequently, $\bigset{\itnwo}_{n\geq K,\,x\in S}$ remains complete. Once again, the removal of an element from a frame either leaves an incomplete set or a frame for $H$ (for example, see \cite[Theorem 8.22]{Hei11}), it then follows that $\bigset{\itnwo}_{n\geq K,\,x\in S}$ is a frame for $H.$ Fix $1\leq i\leq N$ and let $0<\delta<r_i$ be small enough. By Lemma \ref{fraction_conv_operator_norm}, there exists some $K_\delta\in\N$ such that $\frac{\norm{A^{n+1}x}}{\norm{A^nx}}\geq r_i-\frac{\delta}{2}$ for all $n\geq K_\delta$ and all $x\in R_i$. Now let $C$ and $D$ be frame bounds of $\bigset{\frac{A^nx}{\norm{A^nx}}}_{n\geq K_\delta,\,x\in R_i\setminus R_{i-1}}.$  Then for any $y\in H$ such that $\text{supp}(Uy)\subseteq \bigparen{\overline{B_{r_{i}-\delta}}(0)\setminus B_{r_{i-1}+\delta}(0)}$ and any $m\in \N$ we have
\begin{align}
    \begin{split}
        D\,\norm{y}^2&\Ge\sum_{x\in R_i}\sum_{n=K}^\infty\,\Bigabs{\Bigip{y}{\frac{A^nx}{\norm{A^nx}}}}^2 \\
        &\Eq \sum_{x\in R_i}\sum_{n=K}^\infty\,\Bigabs{\Bigip{y}{\frac{A^nx}{\norm{A^nx}}}}^2\\
        &\Eq \sum_{x\in R_i}\sum_{n=K}^\infty\,\Bigabs{\Bigip{(A^*)^{-m}y}{\frac{A^{m+n}x}{\norm{A^nx}}}}^2\quad \text{(By Lemma \ref{invertible_operatorA})}\\
        &\Eq \sum_{x\in R_i}\sum_{n=K}^\infty\,\Bigabs{\Bigip{(A^*)^{-m}y}{\frac{A^{m+n}x}{\norm{A^{m+n}x}}}\frac{\norm{A^{m+n}x}}{\norm{A^{n}x}}}^2
    \end{split}
\end{align}
Note that, by Lemma \ref{limit_operator_norm}, we have that, for any $n\geq K$, $$\frac{\norm{A^{m+n}x}}{\norm{A^nx}}=\frac{\norm{A^{m+n}x}}{\norm{A^{m+n-1}x}}\cdot\frac{\norm{A^{m+n-1}x}}{\norm{A^{m+n-2}x}}\cdots\frac{\norm{A^{n+1}x}}{\norm{A^{n}x}}\geq \Bigparen{\frac{\norm{A^{K+1}x}}{\norm{A^{K}x}}}^m.$$
Consequently, we obtain
\begin{align}
\label{estimate_I}
    \begin{split}
        D\,\norm{y}^2&\Ge \bigparen{r_i-\frac{\delta_1}{2}}^{2m}\sum_{x\in R_i}\sum_{n=K}^\infty\,\Bigabs{\Bigip{(A^*)^{-m}y}{\frac{A^{m+n}x}{\norm{A^{m+n}x}}}}^2\\
        & \Ge \bigparen{r_i-\frac{\delta_1}{2}}^{2m}\Bigparen{C\,\norm{(A^*)^{-m}y}^2-\sum_{x\in R_i}\sum_{n=K}^{m+K-1}\,\Bigabs{\Bigip{(A^*)^{-m}y}{\frac{A^{m+n}x}{\norm{A^{m+n}x}}}}^2}
    \end{split}
\end{align}
Suppose that the closed subspace $M=\clspan{\set{x\in H\,|\,\text{supp}(Ux)\subseteq \overline{B_{r_i-\delta}}(0)\setminus B_{r_{i-1}+\delta}(0)}}$ of $H$ is infinite-dimensional. Then there exists some nonzero $y_0\in M$ such that 
$$\ip{y_0}{A^nx}=0 \quad \text{for all } K\leq n\leq m+K-1 \text{ and all }x\in M.   $$
Moreover, for each $y\in M$ we have $$\norm{(A^*)^{-m}y}^2\Eq \norm{N_\mu^{-m}(Uy)}^2_{\Lm}\Ge \bigparen{r_i-\delta}^{-2m}\norm{y}^2.$$
Thus, it follows from Equation (\ref{estimate_I}) that
$$D\,\norm{y_0}^2\Ge C\bigparen{r_i-\frac{\delta_1}{2}}^{2m}\bigparen{r_i-\delta_1}^{-2m}\norm{y_0}^2_H.$$
Since $m\in\N$ is arbitrary, we obtain a contradiction. Thus, for any $1\leq i\leq N$ and any sufficiently small $0<\delta<r_i$ there exist at most finitely spectral values of $A$ contained in $\overline{B_{r_i-\delta}}(0)\setminus B_{r_{i-1}+\delta}(0)$, which concludes the statements.
%It remains to show that any spectral value of $A$ lying on $\partial B_{r_i}(0)$ is not an eigenvalue of $A$ for all $1\leq i\leq N$. Suppose to the contrary that such eigenvalue $\lambda_i$ exists and let $e_i$ be an eigenvector associated with $\lambda_i.$  Since $\frac{|\lambda_i|^n}{\norm{A^nx}}\geq 1$ for all $x\in S_i$, it follows that 
%$$\sum_{x\in S_i}\sumli\,\Bigabs{\Bigip{e_i}{\frac{A^nx}{\norm{A^nx}}}}^2\Ge |\ip{e_i}{x}|$$
%The statement then follows from that $\sigma(A)\cap B_{\norm{A}}(0)$ contains at most finitely many points for all $0<\delta<\norm{A}$.
\end{proof}

\begin{remark}
Given the assumptions of Theorem \ref{countable_spec_inside}, we note that if the subset $S$ associated with $\mathcal{A}\bigparen{A,S,J_x,(c_{n,x})}$ is such that $\mu\bigparen{\text{supp}(Ux)\cap \partial B_{r_{i}}(0)}=0$ for any $x\notin S_i$, then no spectral value lying on those arcs is an eigenvalue of $A.$ Suppose to the contrary that such eigenvalue $\lambda_i$ exists and let $e_i$ be an eigenvector associated with $\lambda_i.$  By Lemma \ref{complete_full_measure} and the fact that $\frac{|\lambda_i|^n}{\norm{A^nx}}\geq 1$ for all $x\in S_i$, it follows that there must exist some $x\in S_i$ such that
$$\infty>\sumli\,\Bigabs{\Bigip{e_i}{\frac{A^nx}{\norm{A^nx}}}}^2\Ge |\ip{e_i}{x}|^2\sumli|\lambda_i|^{2n}\norm{A^nx}^{-2n}=\infty,$$
which is a contradiction.
\end{remark}

Next, we prove that, under the assumption of Theorem \ref{countable_spec_inside},  $\sigma(A)\cap \partial B_{\norm{A}}(0)$ must have $\mu$-measure zero. Consequently, $A$ must be a diagonal operator if $S$ is a singleton.

\begin{theorem}
\label{diagonal_operator_form}
 Let $A\colon H\rightarrow H$ be a normal operator such that $\mathcal{A}\bigparen{A,S,J_x,(c_{n,x})}$ is a frame for $H$ for some finite subset $S$. 
Assume that there exist some $\eta\in\N$ and $\delta>0$ such that for each infinite $J_x$ there is an increasing syndetic subsequence $(n^x_{k})_{k\in \N}\subseteq J_x$ satisfying $|c_{n^x_{k},x}|\norm{A^{i^x_{k}}x}\geq \delta$ for some non-negative integers $i^x_{k}$ with $|i^x_{k}-n^x_{k}|\leq\eta$ for all $k\in \N$. Then $\mu\bigparen{\sigma(A)\cap \partial B_{\norm{A}}(0)}=0$.

If $|S|=1$, then $A$ is a diagonal operator with all eigenvalues concentrated in $B_{\norm{A}}(0).$

    %Let $A\colon H\rightarrow H$ be a normal operator. Assume that there exists some $x\in H$ such that $\itn$ is a frame for $H$. Then there exists a sequence of distinct complex numbers $\set{\lambda_j}_{j\in \N} \subseteq B_{\norm{A}(0)}$ such that 
    %\begin{enumerate}
     %   \item [\textup{(a)}] $|\lambda_j|\rightarrow  \norm{A}\quad\text{as }j\rightarrow \infty,$
      %  \item [\textup{(b)}] $A=\sum_{j=1}^\infty \lambda_jP_j,$ where $P_j$ are rank one orthogonal projections onto the eigenspace corresponding to $\lambda_j$.
    %\end{enumerate}

\end{theorem}
\begin{proof}
    %Note that by considering $A'=\frac{A}{\norm{A}}$, we may, without loss of generality, assume that $\norm{A}=1$. 
    Likewise, we may assume that $S\cap \text{ker}(A)=\emptyset$. %By Theorem \ref{syndetic_frame_property}, we have that  $\mathcal{A}\bigparen{A,S,\N,(\norm{A^nx}^{-1})}$ is a frame for $H.$
   Suppose that $\text{supp}(\mu)\cap \partial B_{\norm{A}}(0)$ has a positive $\mu$-measure. We then define the closed subspace $$M=\set{x\in H\,|\,Ux=0 ~\mu\text{-a.e. on }B_{\norm{A}}(0)}$$ and let $S_1=\set{x\in S\,|\,(Ux)|_{\partial B_{\norm{A}}(0)}\neq 0\,}$.
    %Without loss of generality, we may assume that $M\neq\set{0}.$ 
    For each $x\in S_1$ we let $x_M$ be the orthogonal projection of $x$ onto $M$. We claim that $\bigset{\frac{(A_M)^nx_M}{\norm{(A_M)^nx_M}}}_{n\geq 1,\,x\in S_1}$ is a frame for $M$, where $A_M$ denote the restriction of $A$ to $M$. Note that $A_M$ is still a normal operator. Clearly, for any $y\in M$ we have 
    
    \begin{equation}
    \label{frame_scaling_ineq}
    \sum_{x\in S}\sum_{n=1}^\infty\, \Bigabs{\Bigip{y}{\frac{(A_M)^nx_M}{\norm{(A_M)^nx_M}}}}^2\Eq \sum_{x\in S_1}\sum_{n=0}^\infty\, \Bigabs{\Bigip{y}{\frac{A^nx\norm{A^nx}}{\norm{(A_M)^nx_M}\norm{A^nx}}}}^2
    \end{equation}
    %\Eq \sum_{n=0}^\infty\, \Bigabs{\Bigip{y}{\frac{A^nx}{\norm{A^nx}}}\frac{\norm{A^nx}}{\norm{(A_S)^nx_S}}}^2.$$
    So, it remains to show that there exist some positive numbers $C$ and $D$ such that $$C\Le \frac{\norm{A^nx}^2}{\norm{(A_M)^nx_M}^2}\Le D \qquad \text{for all } n\geq 0 \text{ and all }x\in S_1.$$
    Clearly, $\frac{\norm{A^nx}^2}{\norm{(A_M)^nx_M}^2}\geq 1$ for all $n\geq 1.$ On the other hand, by Theorem \ref{spectral_thm} and dominate convergence theorem, we have  
    $$\frac{\norm{(A_{M^\perp})^nx_{M^{\perp}}}^2}{\norm{A}^{2n}}\rightarrow 0\qquad \text{as }n\rightarrow\infty$$
    for all $x\in S_1$. It then follows that $$\frac{\norm{A^nx}^2}{\norm{(A_M)^nx_M}^2}=\frac{\norm{(A_M)^nx_M}^2+\norm{(A_{M^{\perp}})^nx_{M^\perp}}^2}{\norm{(A_M)^nx_M}^2}\Eq 1+ \frac{\norm{(A_{M^{\perp}})^nx_{M^\perp}}^2}{\norm{A}^{2n}\norm{x_M}^2}\leq 2$$ for $n$ large enough. By Equation (\ref{frame_scaling_ineq}), we see that $\bigset{\frac{(A_M)^nx_M}{\norm{(A_M)^nx_M}}}_{n\geq 1,\,x\in S_1}$ is a frame for $S$. %Note that $\norm{(A_M)^nx_M}=\norm{A}^n\norm{x_M}$. 
    
    Next, let $A'=\frac{A_M}{\norm{A}}$ and note that 
    $\frac{(A_M)^nx_M}{\norm{(A_M)^nx_M}}=\frac{(A_M)^nx_M}{\norm{A}^{n}\norm{x_M}}$. It then follows that $$\Bigset{(A')^n\Bigparen{\frac{x_M}{\norm{x_M}}}}_{n\geq 0,\,x\in S_1}\Eq \Bigset{\frac{(A_M)^nx_M}{\norm{(A_M)^nx_M}}}_{n\geq 1,\,x\in S_1}$$ is a frame for $M$. However, since $A'$ is a normal operator with $\sigma(A')\subseteq\partial B_{1}(0)$, by Theorem \ref{required_results} (b), this would imply that $M$ is finite-dimensional. Nevertheless, 
     if $M$ were to be a nontrivial finite-dimensional closed subspace of $H$, then there exists at least one eigenvalue $\lambda$ of $A$ with $|\lambda|=\norm{A}$. Let $e_\lambda\in M$ be a nonzero element and note that $\ip{e_\lambda}{x}\neq 0$ for at least one $x\in S_1$ by Lemma \ref{complete_full_measure}. If then follows that
     \begin{align}
         \begin{split}
             \infty > \sum_{x\in S_1}\sum_{n=1}^\infty \,\Bigabs{\Bigip{e_\lambda}{\itnwo}}^2\Eq \sum_{x\in S_1}\sum_{n=1}^\infty \,\frac{\norm{A}^{2n}}{\norm{A^{n}x}^2}\bigabs{\ip{e_\lambda}{x}}^2 \Eq \infty,
         \end{split}
     \end{align}
     which is a contradiction.

     Finally, by Theorem \ref{countable_spec_inside}, we see that $A$ must be a diagonal operator when $S$ is a singleton.
     %Statement (a) now simply follows from the fact that the spectral radius of $A$ equals to $\norm{A}.$
    % \textcolor{red}{Proof need to be changed, the spectrums on $\norm{A}$ may not be point spectrum.}
\end{proof}
\begin{remark} Given the assumptions of Theorem \ref{diagonal_operator_form}, one can show that the restriction of $\mu$ to $\partial B_{\norm{A}}(0)$ is absolutely continuous with respect to the (Lebesgue) arc length measure. By considering $A'=\frac{A}{\norm{A}}$, we may assume that $\norm{A}=1.$ Therefore, for any $y\in H$ we have $$\frac{1}{\norm{x}^2}\sum_{n=1}^\infty\, \bigabs{\bigip{y}{A^nx}}^2\Le \sum_{n=1}^\infty\, \Bigabs{\Bigip{y}{\itnwo}}^2,$$
    which shows that $\set{A^nx}_{n\geq 1}$ is a complete Bessel sequence. The claimed statement then follows from \cite[Theorem 5.1]{ACMCP17}.   \qeddef
\end{remark}
We now apply our main results to partially answer the conjecture proposed by Aldroubi et al.\ that $\mathcal{A}(A,S,\N\cup\set{0},(\norm{A^nx}^{-1}))$ is never a frame for $H$. It is worth noting that this conjecture is not true if $A$ is just a generic bounded linear operator. For example, let $H=\ell^2(\N)$ and let $A\colon H\rightarrow H$ be the right shift operator define by $$A(x_1,x_2,\dots,x_n,\dots)=A(0,x_1,\dots,x_n,\dots)$$
Then $\set{A^n(1,0,0,\dots,)}_{n\geq 0}$ is an orthonormal basis for $\ell^2(\N).$ Moreover, this conjecture also fails if one allows negative powers of $A$ in the iterative system. Let $A\colon L^2([0,1])\rightarrow L^2([0,1])$ be the multiplication operator defined by $(Af)(x)=e^{2\pi ix}f(x)$. Then $A$ is an unitary operator. However, the iterative system $\set{A^nf}_{n\in \Z}$ is a frame for $L^2([0,1])$ for any $f$ such that $|f|$ is bounded above and below by some nonzero constants (for example, see Theorem \cite[Theorem 10.10]{Hei11}).

We say that a sequence of scalars $(c_k)_{k\in \N}$ is \emph{bounded} if there exists some $L>0$ such that $|c_k|\leq L$ for all $k\in \N.$ Applying Theorem \ref{countable_spec_inside} together with Theorem \ref{diagonal_operator_form}, we obtain the following corollary. Note that we proved in proof of Lemma \ref{syndetic_frame_property} that $\Ac\bigparen{A,S,J_x, (c_{n,x})}$ is a frame for $H$ if and only if $\Ac\bigparen{A,S,J_x\cup\set{0}, (c_{n,x})}$ is a frame for $H.$
\begin{corollary} 
\label{partial_solution_Akram_conjecture}
 Let $A\colon H\rightarrow H$ be a normal operator and let $S$ be a finite subset of $H.$ Then for any syndetic subsequences $(n^x_{k})_{k\in\N,\,x\in S}\subseteq \N\cup \set{0}$ and any bounded sequences $(i^x_k)_{k\in\N,\,x\in S}\subseteq \Z$ the iterative system $\Bigset{\frac{A^{n^x_{k}}x}{\norm{A^{n_k^x+i_k^x}x}}}_{k\in\N,\,x\in S}$ is not a frame for $H$ if one of the following three conditions hold:
 \begin{enumerate} 
 \setlength\itemsep{0.5em}
 \item [\textup{(a)}] $\mu\bigparen{\sigma(A)\cap\partial B_{\norm{A}}(0)}\neq 0$
     \item [\textup{(b)}] There exist more than $|S|-1$ complex numbers within $B_{\norm{A}}(0)$ that are contained in the continuous spectrum of $A$ with distinct moduli,
      \item [\textup{(c)}] $S$ is a singleton but $A$ is not a diagonal operator.
 \end{enumerate}
 \qeddef
\end{corollary}
\begin{remark} Assume that $|S|=1$ and assume that $\itn$ is a frame for $H$. Then by Theorem \ref{required_results} (b), Theorem \ref{required_results} (c) and Theorem \ref{diagonal_operator_form}, either the eigenvalues of $A$, $(\lambda_n)\inN,$ are not uniformly separated in $\mathbb{D}$ or the sequence $\bigparen{\norm{P_nx_1}^2(1-|\lambda_n|^{-2})}\inN$ cannot be simultaneously bounded above or below by positive constants. \qeddef
\end{remark}
Finally, we conclude this paper by presenting some illustrative examples.
\begin{example} (a)
Let $H=H^2(\mathbb{D})$ and let $\phi$ be an arbitrary bounded holomorphic function in $\mathbb{D}$. The multiplication operator $A_\phi\colon H^2(\mathbb{D})\rightarrow H^2(\mathbb{D})$ defined by $Af=\phi f$ is a normal operator. Then for any finite subset $S\subseteq H$ and any syndetic subsequences $\set{n_k^x}_{k\in\N,\,x\in S}$ in $ \N\cup\set{0}$ the iterative system $\mathcal{A}\bigparen{A_{\phi},S,\set{n_k^x}_{k\in\N},(\norm{A^{n_k}x}^{-1})}$ is never a frame for $H.$
If $\phi$ is a constant function, then it is clear that the iterative system $\mathcal{A}\bigparen{A_{\phi},S,\set{n_k^x}_{k\in\N},(\norm{A^{n_k}x}^{-1})}$ does not form a frame for $H$. Applying Corollary \ref{partial_solution_Akram_conjecture}, we see that this system also fails to be a frame when $\phi$ is not a constant function.\medskip

(b)  Let $d\geq 1$ be some integer and let $H=L^2(\R^d)$. For any fixed nonzero $\phi \in L^1(\R^d)$, we define the normal operator $A_\phi\colon L^2(\R^d)\rightarrow L^2(\R^d)$ by $A_\phi(f)=\phi\ast f.$ Note that $\widehat{\phi}$ is a non-constant continuous function. Thus, by Corollary \ref{partial_solution_Akram_conjecture}, $\mathcal{A}\bigparen{A_{\phi},S,\set{n_k^x}_{k\in\N},(\norm{A^{n_k}x}^{-1})}$ does not form a frame for $H$ for any finite subset $S\subseteq H$ and any syndetic subsequences $\set{n_k^x}_{k\in\N,\,x\in S}$ in $ \N\cup\set{0}$.

    %Let $H=L^2(\R^d)$ for some $d>1$ and let $\phi\in L^2(\R^d)$ be a bounded function. Then we define $A\colon L^2(\R^d)\rightarrow L^2(\R^d)$
\end{example}

\end{document}